\documentclass{amsart}

\input epsf.tex

\newtheorem{theorem}{Theorem}[section]
\newtheorem{lemma}[theorem]{Lemma}
\newtheorem{prop}[theorem]{Proposition}

\theoremstyle{definition}
\newtheorem{definition}[theorem]{Definition}
\newtheorem{example}[theorem]{Example}

\newtheorem{conjecture}[theorem]{Conjecture}

\theoremstyle{remark}
\newtheorem{remark}[theorem]{Remark}

\numberwithin{equation}{section}

\newcommand{\abs}[1]{\lvert#1\rvert}

\begin{document}
\author[B. Kr\"otz]{Bernhard Kr\"otz}
\author[M. Otto]{Michael Otto}
\address{The Ohio State University, Department of Mathematics, 231 West 18th Avenue, 
Columbus, OH 43210-1174}
\email{kroetz@math.ohio-state.edu, otto@math.ohio-state.edu}
\thanks {The first author was supported in part by NSF-grant DMS-0097314}
\title[Double coset decomposition]
{A Convexity property for the $SO(2,{\mathbb C})$-double coset 
decomposition of $SL(2,{\mathbb C})$ and applications to spherical functions}

\begin{abstract} 
We make a fine study of the $SO(2,{\mathbb C})$-double coset decomposition in $SL(2,{\mathbb C})$
and give a full description of the intersection of the various cells with the complex 
crown $\Xi$ of $SL(2,{\mathbb R})/ SO(2)$. A non-linear convexity theorem is proved 
and applications to analytically continued spherical functions are given. 
\end {abstract}
\maketitle
\section{Introduction}

With a Riemannian symmetric space $X=G/K$ of the non-compact type
comes a natural complexification $\Xi$, the so-called 
{\it complex crown} of $X$. 
One can define $\Xi$ in various ways. To begin 
with let $G=NAK$ be an Iwasawa decomposition and let 
$\mathfrak{a}=\mbox{Lie}(A)$ be the Lie algebra of $A$. 
Studying proper $G$-actions on $X_{\mathbb C}=G_{\mathbb C}/ K_{\mathbb C}$
Akhiezer and Gindikin were led to the following definition 
(cf.\ \cite{AG})  

\begin{equation}\label{a} \Xi=G\exp(i\Omega)K_{\mathbb C}/ K_{\mathbb C}\ ,\end{equation}
where $\Omega\subseteq\mathfrak{a}$ is a certain bounded 
convex set (cf.\ Section 2). Observe that  $\Xi$ is a $G$-invariant 
domain in $X_{\mathbb C}$ containing $X$. Moreover,
the definition of $\Xi$ is independent from the choice of $\mathfrak{a}$
and hence $\Xi$ is generically defined through $X$. 
\par Equivalently, one can define $\Xi$ by 
\begin{equation}\label{b} 
\Xi=\left(\bigcap_{g\in G} g N_{\mathbb C} A_{\mathbb C} K_{\mathbb C}/ K_{\mathbb C} \right)_0\ , 
\end{equation}
where $(\cdot)_0$ refers to the connected component of $(\cdot)$
containing $X$ (cf.\ \cite{FK, GM1, KS1}, \cite{BA, HU},  \cite{FH, GM2, HW}
and \cite{KS3, MA2}).  
Using the {\it complex convexity theorem} from \cite{GK2} one can 
significantly improve on \ref{b}, namely 
\begin{equation}\label{c} 
\Xi=\left(\bigcap_{g\in G} g N_{\mathbb C} A\exp(i\Omega)  K_{\mathbb C}/ K_{\mathbb C} \right)_0 \ . 
\end{equation}
\par The various definitions of $\Xi$ are useful. In \cite{KS2}, where the foundations 
for the $G$-invariant complex geometry of $\Xi$ were laid out 
(construction of natural  $G$-invariant psh exhaustion functions, K\"ahler metrics etc.), one mainly used the 
description (\ref{a}). For other purposes such as analytic continuation of eigenfunctions or the
heat kernel on $X$ to holomorphic  functions on $\Xi$ (cf.\ \cite{KS1}, \cite {KS2})
the characterizations in (\ref{b}) and in particular (\ref{c}) are more appropriate. 

\smallskip One of the main objectives in the study of complex crowns is to achieve 
a better understanding of harmonic analyis on symmetric spaces. For example 
it was shown in \cite{GK1} that all non-compactly causal symmetric 
spaces $G/H$ associated to $G$ appear in the so-called {\it distinguished boundary} of $\Xi$. 
Subsequently this was used in \cite{GKO, GKO2} to construct a Hardy-space for 
non-compactly causal symmetric spaces $G/H$ yielding first progress towards a geometric realization of some of the continuous series in $L^2(G/H)$ (Gelfand-Gindikin program).  

\medskip  The motivation for this article stems from our interest to obtain 
a first understanding of the various kinds of weighted Bergman-spaces which 
can be associated with $\Xi$, most notably Fock-spaces which play 
a prominent role in the context of the heat kernel transform on $X$ (cf.\ \cite{KS2}). 
Intimately related to this circle of problems is the growth behavior of 
analytically continued spherical functions. 
\par A spherical function $\phi_\lambda$
on $X$ is left $K$-invariant, hence its analytic continuation $\tilde\phi_\lambda$  to $\Xi$ 
is necessarily (locally) $K_{\mathbb C}$-invariant. 
This suggests that one should study the relation of $\Xi$ with the complexified polar decomposition 
$\Xi_P\:=K_{\mathbb C} A_{\mathbb C} K_{\mathbb C}/ K_{\mathbb C}$ in $X_{\mathbb C}$. 

\par Whereas the complex crown $\Xi$ is 
contained in the complexified Iwasawa decomposition $N_{\mathbb C} A_{\mathbb C} K_{\mathbb C}/ K_{\mathbb C}$ 
(cf.\ (\ref{b})), the same does not hold for the complexified polar decomposition, i.e., 
$\Xi\not\subseteq \Xi_P$. Furthermore, the 
polar domain $\Xi_P$ is not even open in $X_{\mathbb C}$. These are disappointing 
facts. But before giving up, we at least wanted to understand 
what actually goes ``wrong'' in the most important case of $G=SL(2,{\mathbb R})$. 
It turned out that the above mentioned difficulties  can be overcome, and, moreover, 
by settling the technical problems we could make new 
observations regarding the structure theory of complex crowns. 
\par This paper focuses mostly on $G=SL(2,{\mathbb R})$, but  
we made an effort to present  results in a fashion which allows 
immediate generalization 
to all semisimple Lie groups. An approach in full generality will be given elsewhere. 

\par Let us now describe some of our results in more detail. First observe that 
the polar domain $\Xi_P$ contains a Zariski open subset 
of $X_{\mathbb C}$;  in particular it is dense in $X_{\mathbb C}$. Thus almost all elements in $\Xi$ are 
contained in $\Xi_P$. Our main observation is the following 
non-linear convexity theorem:

\medskip\noindent{\bf Theorem 4.1.} {\it Let $G=SL(2, {\mathbb R})$ and $x=g\exp(iY)K_{\mathbb C} \in \Xi\cap \Xi_P$ with $g\in G$ and 
$Y\in \Omega$ (cf.\ (\ref{a})). Then 
$$x\in K_{\mathbb C} A\exp(i{\rm conv}({\mathcal W}Y))K_{\mathbb C}$$
with ${\rm conv}({\mathcal W}Y)$ the convex hull of the Weyl group orbit of $Y$.}

\medskip We conjecture that Theorem 4.1 holds true for all semisimple Lie groups 
and we provide additional  evidence with a discussion of the Lorentz groups $G=SO_e(1,n)$. 
\par Using general  results of \cite{KS1} one obtains from Theorem 4.1 the following estimate:

\medskip\noindent{\bf Theorem 5.2.} {\it  Let $G=SL(2, {\mathbb R})$. Let $\phi_\lambda$ be a positive definite 
spherical function on $X$ and $\tilde\phi_\lambda$ its analytic continuation 
to $\Xi$. Then $\tilde\phi_\lambda$ is bounded.}

\bigskip Finally,  as an application of our methods we give an 
estimate for the $A_{\mathbb C}$-projection in $\Xi_P$ for 
elements in $\Xi$. 
\par It is our pleasure to thank the referee for his very careful
reading and his useful suggestions. In particular the easy proof
of Theorem \ref{union} is due to him.

\section{Notation and general facts} 

Let $\mathfrak{g} $ be a real semisimple Lie algebra and 
$\ \mathfrak{g}=\mathfrak{k}+\mathfrak{p} \ $ a Cartan 
decomposition. For a maximal abelian subspace $\mathfrak{a}$
of $\mathfrak{p}$ let $\ \Sigma = \Sigma (\mathfrak{g},
\mathfrak{a} ) \subseteq \mathfrak{a}^{*} \ $ be the 
corresponding restricted root system. 
Then $\mathfrak{g}$ admits a root space decomposition
$$ \mathfrak{g} = \mathfrak{a} \oplus \mathfrak{m} \oplus
\bigoplus_{\alpha \in \Sigma} \mathfrak{g}^{\alpha} ,$$
where $\ \mathfrak{m}=\mathfrak{z}_{\mathfrak{k}} (\mathfrak{a})
\ $ and $\ \mathfrak{g}^{\alpha}=\{ X \in \mathfrak{g} :
(\forall H\in \mathfrak{a}) \  \left[ H,X \right] = \alpha(H)X\} \ $ . \\
For a fixed positive system $\Sigma^{+}$ define $\
\mathfrak{n} := \bigoplus_{\alpha \in \Sigma^{+}}
\mathfrak{g}^{\alpha} \ $ . Then we have the Iwasawa 
decomposition on the Lie algebra level:
\[ \mathfrak{g} = \mathfrak{k} \oplus \mathfrak{a} \oplus
\mathfrak{n} . \]
We write $\ \mathcal{ W}=N_{K}(\mathfrak{a})/Z_{K}(\mathfrak{a})
\ $ for the corresponding {\it Weyl group}. 
For an element $\ X \in \mathfrak{a} \ $ we denote by 
$ {\rm conv}(\mathcal{W}X) $ the convex hull of the Weyl group
orbit of $X$. \\
\\
For any real Lie algebra $\mathfrak{l}$ we write 
$\mathfrak{l}_{\mathbb C}$ for its complexification. \\
\par In the sequel $G_{\mathbb C}$ will denote  a simply connected Lie group with Lie
algebra $\mathfrak{g}_{\mathbb C}$. We write  $\ G, K, K_{\mathbb C}, 
A, A_{\mathbb C}, N \ $ and $N_{\mathbb C}$ the analytic 
subgroups of $G_{\mathbb C}$ corresponding to subalgebras $\ \mathfrak{g}, 
\mathfrak{k}, \mathfrak{k}_{\mathbb C}, \mathfrak{a}, 
\mathfrak{a}_{\mathbb C}, \mathfrak{n} \ $ and 
$\mathfrak{n}_{\mathbb C}$ , respectively. \\ 
\\
The following bounded and convex subset of $\mathfrak{a}$ plays a central role:
\[ \Omega := \{ X\in \mathfrak{a} : \abs{\alpha(X)}<
\frac{\pi}{2} \quad \forall \ \alpha \in \Sigma \} . \]
\\
With $\Omega$ we define a left $G$ and right $K_{\mathbb C}$-invariant domain in 
$G_{\mathbb C}$ by 
$$ \tilde \Xi=G\exp(i\Omega)K_{\mathbb C}\ .$$
Also we write 
$$\Xi=\tilde \Xi/ K_{\mathbb C}$$
for the union of right $K_{\mathbb C}$-cosets of $\tilde \Xi$ in the complex symmetric space 
$G_{\mathbb C}/ K_{\mathbb C}$. We refer to $\Xi$ as the {\it complex crown}
of the symmetric space $G/K$. Notice that $\Xi$ is independent of the 
choice of $\mathfrak a$ and hence generically defined through $G/K$.

\par As every root $\alpha \in \Sigma $ 
is analytically integral, it gives rise to a character 
of $A_{\mathbb C}$ by 
$$\ \xi_{\alpha} : A_{\mathbb C} \longrightarrow
\mathbb{C}^* , \ \xi_{\alpha}(\exp(X))=e^{\alpha(X)} \ .$$ \\
Let us define the set of {\it regular elements} in $A_{\mathbb C}$ by 
$$A_{\mathbb{C},\mbox{reg}} = \{ z \in 
A_{\mathbb C} : \xi_{\alpha}^2(z) \neq 1 \quad \forall \ \alpha
\in \Sigma \} .$$
Notice that $A_{\mathbb{C},\mbox{reg}}$ is an algebraic variety
( principal open set in $A_{\mathbb C}$ ). \\

We also define 
$$A_{{\mathbb C}, \mbox{sing}}=A_{\mathbb C}\backslash
A_{\mathbb{C},\mbox{reg}}$$
and call it the {\it singular set} in $A_{\mathbb C}$.

\begin{lemma}
The following assertions hold:
\begin{enumerate}
\item The multiplication mapping
\[ m:K_{\mathbb C} \times A_{\mathbb{C},\mbox{reg}} \times
K_{\mathbb C} \longrightarrow G_{\mathbb C} , \ 
(k_{1},a,k_{2}) \mapsto k_{1}ak_{2} , \]
is submersive. In particular, $\ K_{\mathbb C} 
A_{\mathbb{C},\mbox{reg}} K_{\mathbb C} \subseteq G_{\mathbb C}
\ $ is open.
\item $K_{\mathbb C} A_{\mathbb{C},\mbox{reg}} K_{\mathbb C}$
is dense in $G_{\mathbb C}$.
\end{enumerate}
\end{lemma}
\begin{proof}
(1) is a standard computation which will not be repeated here. 
\par\noindent (2) Notice that $\ K_{\mathbb C} A_{{\mathbb C},\mbox{reg}} K_{\mathbb C}
\subseteq G_{\mathbb C} \ $ is a constructible set as the 
image under the regular mapping 
$\ m:K_{\mathbb C} \times A_{{\mathbb C},\mbox{reg}} \times K_{\mathbb C} 
\longrightarrow G_{\mathbb C} \ $ . Thus
it follows from (1) that $ K_{\mathbb C} A_{{\mathbb C},\mbox{reg}} 
K_{\mathbb C} $ contains a Zariski-open subset of 
$G_{\mathbb C}$. This proves (2). 
\end{proof}

\begin{remark} Notice that the $K_{\mathbb C}$-bi-invariant  domain 
$K_{\mathbb C}A_{\mathbb C} K_{\mathbb C}$ is not open in 
$G_{\mathbb C}$ (cf. our discussion of $G=SL(2,{\mathbb R})$ in Section 3).
\end{remark}

\section{The structure of the double $K_{\mathbb C}$-cosets}
In this section we will give a detailed analysis of the double $K_{\mathbb C}$-coset decomposition 
for $G_{\mathbb C}=SL(2,\mathbb C)$. Further we will start the investigation of the intersections
of $\tilde \Xi$ with the various cells in the double $K_{\mathbb C}$-coset 
decomposition.

\par Let us introduce the necessary notation. Our choice of a maximal compact subgroup of $G$ is $K=SO(2,{\mathbb R})$ and our choice of 
$\mathfrak a$ will be 
$${\mathfrak a}=\left\{ \left( \begin{array}{cc} t & 0 \\ 0 & -t \end{array} \right):\ t\in {\mathbb R}\right\} \ .$$
We choose $\Sigma^+$ such that 
$${\mathfrak n}=\left\{ \left( \begin{array}{cc} 0 & x \\ 0 & 0\end{array} \right):\ x\in {\mathbb R}\right\} \ .$$
Notice that $G_{\mathbb C}=SL(2,{\mathbb C})$ is simply connected and $K_{\mathbb C}=SO(2,{\mathbb C})$.

\par In order to study $K_{\mathbb C}$-double cosets the use of spherical functions is 
useful. 

\begin{definition} On $G_{\mathbb C}=SL(2,{\mathbb C})$ we define the {\it elementary spherical function} $\Phi$
by 
\[ \Phi :G_{\mathbb{C}} \longrightarrow \mathbb{C}, \quad 
\Phi (g) := \mbox{tr}(gg^{t}) \quad \forall \ g\in G_{\mathbb{C}} . \]
( For $\ g = \left( \begin{array}{cc} a & b \\ c & d \end{array} \right)\ $ we have $\ \Phi(g)=a^2 +b^2 +c^2+d^2$. )
\par Note that the function $\Phi $ is holomorphic and 
$K_{\mathbb C}$-bi-invariant.
\end{definition}
\par The set of singular elements in $A_{\mathbb C}$ is given by 

$$A_{{\mathbb C}, \mbox{sing}}=\left\{ \left( \begin{array}{cc} z & 0 \\ 0 & z^{-1} \end{array} \right) \
:\ z\in\{-i,i,-1,1\}\right\} \ .$$

\par Let $\ g \in G_{\mathbb C} \setminus K_{\mathbb C} 
A_{\mathbb C} K_{\mathbb C} . \ $ According to Lemma 2.1 
we can find sequences \[
(k_{n}) \subset K_{\mathbb C}, \ (a_{n}) \subset A_{\mathbb C}
, \ (\tilde{k}_{n}) \subset K_{\mathbb C} \quad
\mbox{with} \quad k_{n}a_{n}\tilde{k}_{n} \longrightarrow g . 
\]

\begin{prop}\label{limits}
Let $\ g \in G_{\mathbb C} \setminus K_{\mathbb C} A_{\mathbb C} K_{\mathbb C} \ $ be the limit of the sequence $\ 
(k_{n}a_{n}\tilde{k}_{n})
\subset K_{\mathbb C} A_{\mathbb C} K_{\mathbb C} . $ 
Then  $(a_{n})\subset A_{\mathbb C}$ is bounded and for every 
convergent subsequence $(a_{n_k})$ we have  
$$\lim_{k \to \infty} a_{n_k} \in A_{{\mathbb C}, \mbox{sing}}\ . $$
\end{prop}
\begin{proof}
Write $\ a_{n}=\left( \begin{array}{cc} z_{n} & 0 \\ 0 & z_{n}^{-1} \end{array} \right) \ $. 
By continuity of the map $\Phi $ we get
\[ \Phi (g) = \mbox{tr}(gg^{t}) = \lim_{n \to \infty } \mbox{tr}(a_{n}a_{n}^{t})
= \lim_{n \to \infty} (z_{n}^{2}+z_{n}^{-2}) . \]
It follows that there exist $\ m, M > 0 \ $ such that
\[ 0 < m < \abs{z_{n}} < M < \infty \quad \forall n . \]
Passing to an appropriate subsequence we can assume that $\ 
z_{n} \longrightarrow z_{0} \neq 0 . $ \\ 
\\
Suppose $(k_{n})$ has a convergent subsequence. W.l.o.g. we may assume that  \[
k_{n} \longrightarrow k_{0} 
\in K_{\mathbb C}, \quad a_{n} \longrightarrow a_{0} \in A_{\mathbb C} . \]
But then $\ \tilde{k}_{n} = a_{n}^{-1} k_{n}^{-1} (k_{n}a_{n}\tilde{k}_{n}) \
$ also converges. \\
We have now
\[ k_{n} \longrightarrow k_{0} 
\in K_{\mathbb C}, \quad a_{n} \longrightarrow a_{0} \in A_{\mathbb C}, 
\quad \tilde{k}_{n} \longrightarrow \tilde{k}_{0} \in K_{\mathbb C} . \]
But this would imply $\ g = \lim_{n\to\infty} k_{n} a_{n} \tilde{k}_{n} \in 
K_{\mathbb C} A_{\mathbb C} K_{\mathbb C} , \ $ contradicting the assumption.
\\
By the same argument we see that $(\tilde{k}_{n})$ cannot have a 
convergent subsequence. \\ 
\\
We want to know more about the limits $z_{0}$. \\
Notice that $k_n a_n^2 k_n^t$ converges in $G_{\mathbb C}$. 
For simplicity we omit the indices writing $k=k_{n}$ and
$a=a_{n}$. We look at the elements $\ ka^{2}k^{t} \in G_{\mathbb C}/
K_{\mathbb C} . $  \\ 
Due to the $K_{\mathbb C}$-bi-invariance of $ K_{\mathbb C} A_{\mathbb C} 
K_{\mathbb C} $ and the compactness of $K$ we can assume that 
\[ k=\left( \begin{array}{cc}
\cosh t & i\sinh t \\ -i\sinh t & \cosh t \end{array} \right) 
= \left( \begin{array}{cc}
\cos it & \sin it \\ -\sin it & \cos it \end{array} \right),  \]
where $t\in {\mathbb R}$ and $|t|\to \infty$. 
Then with $a=\left( \begin{array}{cc}
z& 0 \\ 0 & z^{-1} \end{array} \right)$ we have 
\begin{eqnarray}\label{matrix}
ka^{2}k^{t} & = & \left( \begin{array}{cc}
\cosh t & i\sinh t \\ -i\sinh t & \cosh t \end{array} \right) 
\left( \begin{array}{cc}
z^{2} & 0 \\ 0 & z^{-2} \end{array} \right)
\left( \begin{array}{cc}
\cosh t & -i\sinh t \\ i\sinh t & \cosh t \end{array} \right) 
\nonumber \\
& = &
\left( \begin{array}{cc}
z^{2} \cosh^{2}t-z^{-2}\sinh^{2}t  & (z^{-2}-z^{2})i\sinh t \cosh t \\
(z^{-2}-z^{2})i\sinh t \cosh t & z^{-2} \cosh^{2}t -z^{2} \sinh^{2} t
\end{array} \right) . 
\end{eqnarray}
We know that $\ z \to z_{0} \ $ and $\ \abs{t} \to \infty . $ \\
The upper left entry in (\ref{matrix}) must converge, therefore
\[ \frac{z^{2} \cosh^{2}t-z^{-2}\sinh^{2}t}{\sinh^{2}t} = 
z^{2} \frac{\cosh^{2}t}{\sinh^{2}t}-z^{-2} \longrightarrow 0 . \]
But this implies $\ z_{0}^{2}=z_{0}^{-2} \ $ , hence $\ z_{0}^{2}=\pm 1 . 
$ 
\end{proof}

\par The vector space $\mathfrak{p}_{\mathbb C}$ admits a decomposition 
$\mathfrak{p}_{\mathbb C}=\mathfrak{p}^+\oplus \mathfrak{p}^-$
into irreducible $\mathfrak{k}_{\mathbb C}$-modules with 

\[ \mathfrak{p}^{-} = \mathbb{C} \left( \begin{array}{cc}
1 & -i \\ -i & -1 \end{array} \right) \quad \mbox{and} \quad
\mathfrak{p}^{+} = \mathbb{C} \left( \begin{array}{cc}
1 & i \\ i & -1 \end{array} \right) . \]
The corresponding analytic subgroups of $G_{\mathbb C}$ are given by 
\[ P^{-} = \left\{ \left( \begin{array}{cc} 1-u & iu \\ iu & 1+u
\end{array} \right) : u \in \mathbb C \right\} \quad \mbox{and} \quad
P^{+} = \left\{ \left( \begin{array}{cc} 1+u & iu \\ iu & 1-u
\end{array} \right) : u \in \mathbb C \right\} . \]
We introduce punctured discs in $P^-$ and $P^+$ by 
$$P_{1\over 2}^-=\left\{ \left( \begin{array}{cc} 1-v & iv \\ iv & 1+v \end{array} 
\right) : v \in \mathbb{C}, 0<\abs{v}<\frac{1}{2} \right\} $$
and similarily we define $P_{1\over 2}^+$. \\
Finally, we set
\[ A_{\mathbb C}':= A_{\mathbb C} \setminus \{ \pm {\bf 1} \} \quad
\mbox{and} \quad P^{\mp '}:= P^{\mp} \setminus \{ {\bf 1} \} , \]
\\
and define an element 
$$y_0=\left( \begin{array}{cc} i & 0 \\ 0 & -i
\end{array} \right) \in A_{{\mathbb C},\mbox{sing}}\ .$$
We can now describe how $G_{\mathbb C}$ decomposes into a disjoint
union of $K_{\mathbb C}$-bi-invariant subsets. 
Furthermore, we will give a complete description of the 
intersection of $\tilde\Xi=G\exp(i\Omega)K_{\mathbb C}$ with the lower 
dimensional cells. The 
intersection of $\tilde\Xi$ with the big cell $K_{\mathbb C} A_{\mathbb C} 
K_{\mathbb C} $ will be subject of the next section.
\\
\\
\begin{theorem}\label{union}
\begin{enumerate}
\item The $K_{\mathbb C}$-double coset decomposition of $G_{\mathbb C}$ is given  by:
\[ G_{\mathbb C}=K_{\mathbb C} A_{\mathbb C}' K_{\mathbb C} \
\amalg \ P^{-'} K_{\mathbb{C}} \ \amalg \ 
P^{+'} K_{\mathbb{C}} \ \amalg \ P^{-'} y_0 K_{\mathbb{C}} \ 
\amalg \ P^{+'} y_0 K_{\mathbb{C}} \ \amalg \ K_{\mathbb C}
. \]
\item The following equality holds:
$$\tilde\Xi  \setminus K_{\mathbb C} A_{\mathbb C} 
K_{\mathbb C} = P_{1\over 2}^-K_{\mathbb C} \amalg P_{1\over 2}^+K_{\mathbb C}\ .$$
\end{enumerate}
\end{theorem}
\medskip
In order to prove Theorem \ref{union} it is useful to adopt a more geometric point of view. 
We closely follow the suggestions of the referee. 
\par In the sequel we realize $X=SL(2,\mathbb{R})/SO(2, \mathbb{R})$ as 
the upper sheet of the two-sheeted hyperboloid
\[ x_{0}^{2}-x_{1}^{2}-x_{2}^{2}=1, \ x_{0}>0 \quad (x=(x_0, x_1, x_2)\in \mathbb{R}^{3}) \]
and $X_{\mathbb{C}}=SL(2,\mathbb{C})/SO(2,\mathbb{C})$
as the complex quadric 
\[ z_{0}^{2}-z_{1}^{2}-z_{2}^{2}=1 \quad (z=(z_0,z_1,z_2)\in \mathbb{C}^{3}). \]
Let $e_0=(1,0,0)$. Then the polar decomposition in $X$ is given by $X=KA. e_{0}$,
\[ x_{0}=\cosh t, \ x_{1}=\sinh t \cos \phi , \ x_{2}=\sinh t \sin \phi  \]
with $t,\varphi\in \mathbb{R}$. 
\medskip

\begin{proof}
$(Thm.\ref{union})$
Complexifying the real polar decomposition of $X$ it follows that a point 
$z\in X_{\mathbb{C}}$ belongs to $K_{\mathbb{C}}A_{\mathbb{C}}. e_{0}$
if and only if there exist $a,b,u,v \in \mathbb{C}$ with
\[ a^{2}+b^{2}=1, \quad u^{2}-v^{2}=1 , \]
such that
\[ z_{0}=u, \ z_{1}=va, \ z_{2}=vb . \]
Therefore $X_{\mathbb{C}}\setminus K_{\mathbb{C}}A_{\mathbb{C}}. e_{0} $
decomposes into four $K_{\mathbb{C}}$-orbits
\begin{eqnarray*}
(1) & & z_{0}=1, z_{1}=\tau, z_{2}=i\tau, \\
(2) & & z_{0}=1, z_{1}=\tau, z_{2}=-i\tau, \\
(3) & & z_{0}=-1, z_{1}=\tau, z_{2}=i\tau, \\
(4) & & z_{0}=-1, z_{1}=\tau, z_{2}=-i\tau, 
\end{eqnarray*}
where $\tau \in \mathbb{C}^{*}$. \\
The orbits $(1),(2),(3),(4)$ correspond to the $K_{\mathbb{C}}$-bi-invariant 
sets $P^{-'} K_{\mathbb{C}}$, $P^{+'} K_{\mathbb{C}}$, $P^{-'} y_0 K_{\mathbb{C}}$,
$P^{+'} y_0 K_{\mathbb{C}} $, respectively. \\
This proves part $(1)$ of the theorem. \\
\\
The points $z=x+iy \in \Xi \subset X_{\mathbb{C}}$ are characterized by
the property
\[ x_{0}>0, \quad x_{0}^{2}-x_{1}^{2}-x_{2}^{2}>0 . \]
For an element $z=(1,\tau,i\tau)$ in the orbit $(1)$ this translates to
$\abs{\tau}<1$. Therefore,
\[ \{ (1,\tau,i\tau)\in X_{\mathbb{C}} : \tau \in \mathbb{C}^{*} \}
\cap \Xi = \{ (1,\tau,i\tau)\in X_{\mathbb{C}} : \tau \in \mathbb{C}^{*},
\ \abs{\tau}<1 \} \ , \]
and this subset of $\Xi$ corresponds to $P_{1\over 2}^-K_{\mathbb C} \subset
\tilde{\Xi} $. \\
Similarly one relates orbit $(2)$ with $P_{1\over 2}^{+}K_{\mathbb C} \subset
\tilde{\Xi} $. \\
The intersections of both orbits $(3)$ and $(4)$ with $\tilde \Xi$ are empty.
\end{proof}

\begin{remark} (a) Theorem \ref{union}(1) can be deduced 
from  general results of Matsuki (cf.\ \cite{MA1}).
\par\noindent (b) The cells in the decomposition of $G_{\mathbb C}$ in Theorem \ref{union}
are characterized by the values of $\Phi $. 
The function $\Phi $ attains the value $2$ on $K_{\mathbb C}$ and on $\ P^{-'} 
K_{\mathbb{C}} \amalg P^{+'} K_{\mathbb{C}} \ $, the value $-2$ on $ \
P^{-'} y_0 K_{\mathbb{C}} \amalg P^{+'} y_0 K_{\mathbb{C}} \ $, and $\Phi $
can take any value except $2$ on $K_{\mathbb C}A_{\mathbb C}'K_{\mathbb C}$.
\end{remark}

\section{A non-linear convexity theorem}
In this section we will investigate the intersection of $\tilde\Xi$ with the big cell 
$K_{\mathbb C} A_{\mathbb C} K_{\mathbb C}$. 
An element $\ g \in
 G \exp (i\Omega) K_{\mathbb{C}} \cap K_{\mathbb C} 
A_{\mathbb C} K_{\mathbb C} \ $ can be written as 

\begin{equation}\label{g}
g = h \left( \begin{array}{cc} e^{i\theta} & 0 \\ 0 & 
e^{-i\theta} \end{array} \right) k = \tilde{k_{1}} \tilde{a} 
\tilde{k_{2}} , 
\end{equation}
where $\tilde k_1, \tilde k_2\in K_{\mathbb C}$ and 
\[ h = \left( \begin{array}{cc} a & b \\ c & d
\end{array} \right) \in G \quad \mbox{and} \quad \tilde{a}=\left( 
\begin{array}{cc} z & 0 \\ 0 & z^{-1} \end{array} \right)
\in A_{\mathbb C} . \]

We want to describe the dependence of $\tilde{a}$ from $\theta$ and $h$.
\\
\\
Recall that for $SL(2,\mathbb{R})$ the Weyl group consists of only two
elements and for $X\in \mathfrak{a}$ we have $\ {\rm conv}
(\mathcal{W}.X) = [-1,1] \cdot X \ $.

\begin{theorem}\label{conv}
Let $ X \in \Omega $ and $ h \in G $. Suppose that $\ 
g = h \exp(iX) = \tilde{k_{1}} \tilde{a} \tilde{k_{2}} \in 
K_{\mathbb C} A_{\mathbb C} K_{\mathbb C} \ $. Then
\[ g \in K_{\mathbb C} A \exp(i{\rm conv}(\mathcal{W}X)) 
K_{\mathbb C}
. \]
\end{theorem}
\begin{proof}
Using the Iwasawa decomposition $G=KAN$ we may assume that
\[
\begin{array}{cccccccc}
g & = & \underbrace{\left( \begin{array}{cc} 1 & r \\ 0 & 1 \end{array} 
\right) } & \underbrace{a\exp iX } & = & \underbrace{\left( \begin{array}{cc}
\cos \omega & \sin \omega \\ -\sin \omega & \cos \omega \end{array} 
\right) } & \underbrace{\left( \begin{array}{cc} z & 0 \\ 0 & z^{-1} 
\end{array} \right) } & \underbrace{\tilde{k}_{2} } \\
& & \in N & \in A\exp (i\Omega) & & \in K_{\mathbb C} & 
=\tilde{a} \in A_{\mathbb C} &
\in K_{\mathbb C}  
\end{array} \]
If we set $\ \left( \begin{array}{cc} u & 0 \\ 0 & u^{-1} 
\end{array} \right) := a\exp iX \ $ and $\ \theta := \mbox{Arg}(u) \ $,
then
\[ gg^{t}= \left( \begin{array}{cc}
u^{2}+r^{2}u^{-2} & ru^{-2} \\ ru^{-2} & u^{-2} \end{array} \right) =
\left( \begin{array}{cc} z^{2}\cos^{2}\omega + z^{-2}\sin^{2} \omega & 
(z^{-2}-z^{2})\sin \omega \cos \omega \\ (z^{-2}-z^{2})\sin \omega \cos 
\omega & z^{2}\sin^{2} \omega + z^{-2}\cos^{2} \omega \end{array} \right) ,
\] yielding 
\begin{equation}\label{form1}
\Phi (g)=\mbox{tr}(gg^{t})=u^{2}+(1+r^{2})u^{-2} = z^{2}+z^{-2} . 
\end{equation}
In Figure 1 the set of points of the hyperbola is
$\ \{ l e^{2i\theta} + l^{-1} e^{-2i\theta} : l \in \mathbb{R}, l>0
\} \ $. The complex number $z^{2}+z^{-2}$ lies within the shaded region. 
Therefore, there exists $w=l e^{i\theta}$ such that 
$\ z^{2}+z^{-2} = \lambda (w^{2}+w^{-2}) \ $ for some real number $\lambda >1$.

\medskip
\epsfbox{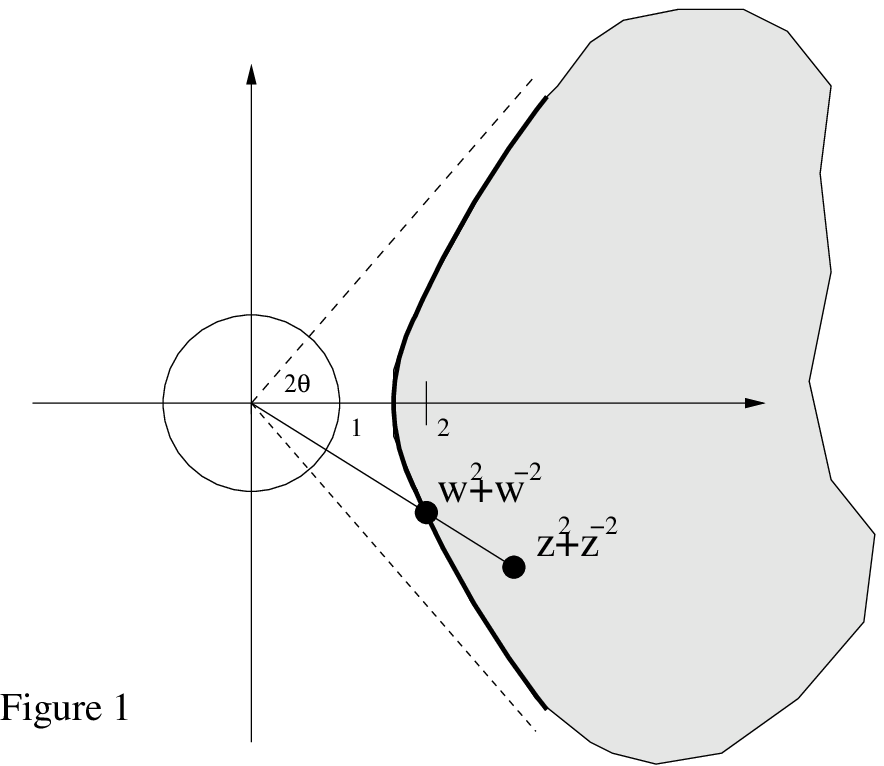}
\medskip 
\par 
We now obtain the assertion of the theorem with the following Lemma \ref{convlem}.  
\end{proof}

\begin{lemma}\label{convlem} Let $r_1, r_2>0$ and $-{\pi\over 2}<\varphi_1,\varphi_2<{\pi\over 2}$.
Define complex numbers by $u=r_{1} e^{i\varphi_{1}}$ and $v=r_{2} e^{i\varphi_{2}}$. 
Assume that $\ v+v^{-1} = \lambda (u+u^{-1}) \ $, with $\lambda \in \mathbb{R}, \lambda
>1 $. Then, $\quad \abs{\varphi_{2}} \leq \abs{\varphi_{1}} \ $.
\end{lemma}
\begin{proof}
Clearly we can assume $\varphi_{2} \neq 0 $. 
Comparing real and imaginary parts of $\ v+v^{-1} = \lambda 
(u+u^{-1}) \ $ , followed by squaring, yields
\[  \quad \lambda^{2}(r_{1}+r_{1}^{-1})^{2}\cos^{2}{\varphi_{1}} = (r_{2}+r_{2}^{-1})^{2}\cos^{2}{\varphi_{2}} \ , \] and
\begin{equation}\label{eq1}
\lambda^{2}(r_{1}-r_{1}^{-1})^{2}\sin^{2}{\varphi_{1}} = (r_{2}-r_{2}^{-1})^{2}\sin^{2}{\varphi_{2}} \ . 
\end{equation} 
Adding these two equations we get 
\begin{equation}\label{eq2} 
\lambda^{2}(r_{1}-r_{1}^{-1})^{2}+4\lambda^{2}\cos^{2}{\varphi_{1}} = 
(r_{2}-r_{2}^{-1})^{2} + 4 \cos^{2}{\varphi_{2}} \ . 
\end{equation} 
Combining (\ref{eq1}) and (\ref{eq2}) gives

\begin{equation}\label{eq3} 
\lambda^{2}(r_{1}-r_{1}^{-1})^{2}+4\lambda^{2}\cos^{2}{\varphi_{1}} = \lambda^{2}(r_{1}-r_{1}^{-1})^{2}
\frac{\sin^{2}{\varphi_{1}}}{\sin^{2}{\varphi_{2}}} + 4\cos^{2}{\varphi_{2}} \ .
\end{equation}
If we assume $\ \abs{\varphi_{2}} > \abs{\varphi_{1}} \ $, then 
\[ \lambda^{2}(r_{1}-r_{1}^{-1})^{2} > \lambda^{2}(r_{1}-r_{1}^{-1})^{2} \frac{\sin^{2}{\varphi_{1}}}{\sin^{2}{\varphi_{2}}}
, \] 
and \[ 4\lambda^{2}\cos^{2}{\varphi_{1}} > 4\cos^{2}{\varphi_{2}} \ , \] 
a contradiction to (\ref{eq3}). 
\end{proof}

\begin{conjecture} Theorem \ref{conv} holds for any real connected semisimple
Lie group $G$. \end{conjecture}

\begin{example} To give more evidence for our conjecture we will now discuss the case 
of $G=SO_{e}(1,n)$, $n\geq 2$. 
\par On the Lie algebra level we have the Cartan decomposition $ \mathfrak{g}=
\mathfrak{k} + \mathfrak{p} $, where
\begin{eqnarray*}
\mathfrak{k} & = & \left\{ \left( \begin{array}{cc} 0 & 0 \\ 0 & B \end{array} 
\right) \in \mbox{Mat}(n+1,\mathbb{R}) : B \in \mathfrak{so}(n) \right\} \ ,  \\ 
\mathfrak{p} & = &\left\{ \left( \begin{array}{cc} 0 & v \\ v^{t} & 0 \end{array} \right)
\in \mbox{Mat}(n+1,\mathbb{R}) : v \in \mathbb{R}^{n} \right\} . 
\end{eqnarray*}
Choosing 
\[ \mathfrak{a} = \left\{ \left( \begin{array}{ccc} 0 & 0 & t \\ 0 & 0 & 0 \\
t & 0 & 0 \end{array} \right) \in \mathfrak{p} : t \in \mathbb {R} \right\}  \]
as maximal subalgebra of $\mathfrak{p}$ yields the root system $\ \Sigma =
\{ \alpha , -\alpha \} \ $ with 
\[ \alpha \left( \begin{array}{ccc} 0 & 0 & 1 \\ 0 & 0 & 0 \\
1 & 0 & 0 \end{array} \right) = 1. \]
Therefore,
\[ \Omega = \{ X\in \mathfrak{a} : \abs{\alpha(X)}<
\frac{\pi}{2} \quad \forall \ \alpha \in \Sigma \} = \left\{ \left( 
\begin{array}{ccc} 0 & 0 & t \\ 0 & 0 & 0 \\ t & 0 & 0 \end{array} \right) :
t\in \mathbb{R}, \abs{t}<\frac{\pi}{2} \right\} . \]
We choose $\alpha $ as the positive root. Then the subgroups in the 
Iwasawa decomposition $G=KAN$ are 
\begin{eqnarray*}
K & = & \left\{ \left( \begin{array}{cc} 1 & 0 \\ 0 & B \end{array} 
\right) : B \in SO(n) \right\} \ \simeq \ SO(n) , \\
N & = & \left\{ \left( \begin{array}{ccc} 1+\frac{1}{2}\abs{\abs{v}}^{2} & v & 
-{1\over 2}\abs{\abs{v}}^{2} \\ v^{t} & \mbox{Id}_{n-1} & -v^{t} \\
\frac{1}{2}\abs{\abs{v}}^{2} & v & 1-\frac{1}{2}\abs{\abs{v}}^{2} 
\end{array} \right) : v \in \mathbb{R}^{n-1} \right\} , \\
A & = & \left\{ \left( \begin{array}{ccc} \cosh t & 0 & \sinh t \\ 0 & 
\mbox{Id}_{n-1} & 0 \\ \sinh t & 0 & \cosh t \end{array} \right) : t \in
\mathbb{R} \right\} .
\end{eqnarray*}
We consider the standard representation $(\pi ,V)$ of $G$ on $V=\mathbb{C}^{
n+1}$. \\
Clearly, $\ v_{0}=e_{1}=(1,0,\cdots ,0)^{t} \ $ is a $K$-fixed vector. The 
elementary spherical function $\Phi $ on $\ G_{\mathbb C}\simeq 
SO(n+1,\mathbb{C}) \ $ is given by
\[ \Phi :G_{\mathbb C} \longrightarrow \mathbb{C}, \quad \Phi (g)= \langle \pi(g)
v_{0},v_{0}\rangle  . \]
Note that $\Phi(g)$ is just the upper left entry $g_{11}$ in the matrix
$(g_{ij})_{1\leq i,j \leq n+1} $. \\
Now let $\ x \in G\exp(i\Omega) \cap K_{\mathbb C}A_{\mathbb C}K_{\mathbb C}
\ $. W.l.o.g. we can write 
\[ x=na\exp iX = \tilde{k}_{1} \tilde{a} \tilde{k}_{2} , \]
for some
\begin{eqnarray*}
n  & = & \left( \begin{array}{ccc} 1+\frac{1}{2}\abs{\abs{v}}^{2} & v & 
-{1\over 2}\abs{\abs{v}}^{2} \\ v^{t} & \mbox{Id}_{n-1} & -v^{t} \\
\frac{1}{2}\abs{\abs{v}}^{2} & v & 1-\frac{1}{2}\abs{\abs{v}}^{2} 
\end{array} \right) \in N ,  \\
a\exp iX   & = &   \left( \begin{array}{ccc} \cosh z & 0 & \sinh z \\ 0 & 
\mbox{Id}_{n-1} & 0 \\ \sinh z & 0 & \cosh z \end{array} \right) \in A
\exp (i\Omega) \\
\end{eqnarray*}
and 
$$\tilde{a}= \left( \begin{array}{ccc} \cosh \omega & 0 & \sinh \omega \\ 0 & 
\mbox{Id}_{n-1} & 0 \\ \sinh \omega & 0 & \cosh \omega \end{array} \right)  
\in A_{\mathbb C} , \quad
\tilde{k}_{1,2}=\left( \begin{array}{cc} 1 & 0 \\ 0 & B_{1,2} \end{array} 
\right) \ ( \ B_{1,2} \in SO(n,\mathbb{C}) \ ) \ .$$
Applying $\Phi $ yields
\[ \Phi(x)=(1+\frac{1}{2}\abs{\abs{v}}^{2}) \cosh z - 
\frac{1}{2}\abs{\abs{v}}^{2} \sinh z = \cosh \omega , \]
or
\begin{equation}\label{SO(1,n)}
e^{z}+ (1+\abs{\abs{v}}^{2})e^{-z} = e^{\omega}+e^{-\omega} . 
\end{equation}
Equation (\ref{SO(1,n)}), Figure 1 (with a slight modification in the 
notation) and Lemma \ref{convlem}  show again that $\ \tilde{a} \in A\exp(i{\rm conv}(\mathcal{W}.X))
K_{\mathbb C} \ $, as claimed. 
\end{example}

\section{Applications to spherical functions}

Let us briefly introduce the necessary notation. Write $G=NAK$ for the Iwasawa
decomposition of $G=SL(2,{\mathbb R})$ with $N$, $A$ and $K$ as in Section 3. 
For $g\in G$ we define $a(g)\in A$ by $g\in Na(g)K$. For $\lambda\in {\mathfrak a}_{\mathbb C}^*$
and $a\in A$ we set $a^\lambda=e^{\lambda(\log a)}$. 
\par Following Harish-Chandra we define the {\it spherical function on $G/K$ with parameter 
$\lambda\in {\mathfrak a}_{\mathbb C}^*$} by 

$$\phi_\lambda(gK)=\int_K a(kg)^{\rho-\lambda}\ dk \qquad (g\in G), $$
where $dk$ denotes the normalized Haar-measure on $K$ and $\rho={1\over 2}\alpha$
with $\Sigma^+=\{\alpha\}$. 
\par Define 
$$\Pi=\{\lambda\in \mathfrak{a}_{\mathbb C}^*:\  \phi_\lambda \ \mbox{is positive definite}\}$$
and recall that $i{\mathfrak a}^*\subseteq \Pi$. 
\par We have to recall some facts from \cite{KS1} \S4  on 
the analytic continuation of the spherical functions (actually 
valid for all semisimple Lie groups $G$).

\begin{prop}\label{KS}
 Let $\lambda\in {\mathfrak a}_{\mathbb C}^*$. Then the following assertions hold:
\begin{enumerate}
\item  The spherical function $\phi_\lambda$ has a unique holomorphic extension 
to $\Xi$. 
\item The spherical function $\phi_\lambda$ has a unique continuation 
to a $K_{\mathbb C}$-bi-invariant function on $K_{\mathbb C} A\exp(2i\Omega)K_{\mathbb C}/ K_{\mathbb C}$
such that the restriction to $A\exp(2i\Omega)$ is holomorphic. 
Moreover, if $\lambda\in \Pi $ and $\Omega_c\subseteq\Omega$ is 
a compact subset, then the restriction of $\phi_\lambda$ to 
$K_{\mathbb C} A\exp(2i\Omega_c)K_{\mathbb C}/ K_{\mathbb C}$ is bounded. 
\end{enumerate}
\end{prop}

In the sequel we also denote by $\phi_\lambda$ the holomorphic extension of $\phi_\lambda$
to $\Xi$. We denote by $Z=\{\pm{\bf 1}\}$ the center of $G$. 
Our result in this section then is:

\begin{theorem} \label{sph} 
Let $\lambda\in \Pi$. Then the spherical function $\phi_\lambda$ is 
bounded on $\Xi$, i.e., 
$$\|\phi_\lambda\|_{\infty, \Xi}<\infty\ .$$
\end{theorem}

\begin{proof} By Theorem \ref{conv} we have 
$$\Xi\bigcap K_{\mathbb C} A_{\mathbb C} K_{\mathbb C}/ K_{\mathbb C} \subseteq K_{\mathbb C} A\exp(i\Omega)K_{\mathbb C}/ K_{\mathbb C}\ .$$
Thus Proposition \ref{KS} implies that $\phi_\lambda$ is bounded 
on $\Xi\bigcap K_{\mathbb C} A_{\mathbb C} K_{\mathbb C}/ K_{\mathbb C}$. 
\par Let now $z\in \Xi$ such that $z\not\in K_{\mathbb C} A_{\mathbb C} K_{\mathbb C}/ K_{\mathbb C}$. 
By the density of $K_{\mathbb C} A_{\mathbb C} K_{\mathbb C}$ in $G_{\mathbb C}$
we can find sequences $(k_n)\subseteq K_{\mathbb C}$, 
$(a_n)\subseteq A_{\mathbb C}$ such that $k_n a_n K_{\mathbb C} \in \Xi$ for all
$n\in {\mathbb N}$ and $z=\lim_{n\to\infty} k_n a_n K_{\mathbb C}$. By Theorem \ref{conv} 
we have $a_n\in ZA\exp(i\Omega)$ and thus 
$$\lim_{n\to \infty} a_n=\pm {\bf 1}$$
by Proposition \ref{limits}.  By the continuity of $\phi_\lambda$ we hence obtain 
that 
$$\phi_\lambda(z)=\lim_{n\to \infty}\phi_\lambda(k_n a_n K_{\mathbb C})
=\lim_{n\to \infty}\phi_\lambda(a_n K_{\mathbb C})=\phi_\lambda(\pm {\bf 1})=1, $$
concluding the proof of the theorem.
\end{proof}

\begin{remark} 
(a) Note that Proposition \ref{KS} holds for all 
semisimple Lie groups. Thus in order to generalize Theorem \ref{sph} 
to all semisimple Lie groups one needs Proposition 
\ref{limits} and Theorem \ref{conv} in full generality (cf. our conjecture 
in Section 4). 

\par\noindent  (b) Fix a $\lambda\in {\mathfrak a}_{\mathbb C}^*$  
and consider the following function on $G/K\times K$
$${\mathcal P}_\lambda(gK, k)\mapsto a(kg)^{\rho-\lambda}\ .$$
Identifying $G/K$ with the unit disc  $D$ and $K$ with its boundary, this 
function is easily seen to be a power of the Poisson kernel on $D$. According 
to \cite{KS1} the function ${\mathcal P}_\lambda$ admits an analytic continuation to 
$\Xi\times K$ which is holomorphic in the first variable. 
Also by \cite{KS1}  we know that this function is unbounded for $\lambda\in \Pi\setminus\{\rho, -\rho\}$. 
Observe 
that the spherical function on $\Xi$ is given by 
$$\phi_\lambda(z)=\int_K {\mathcal P}_\lambda(z, k) \ dk \qquad (z\in \Xi)\ .$$
Now notice that Theorem \ref{sph} says that $\phi_\lambda$ stays 
bounded despite the fact that the integrands become singular 
towards the boundary of $\Xi$.
\end{remark}

\section{A norm estimate for the middle projection}

In this section we will investigate the growth of the 
middle projection in the $K_{\mathbb C} A_{\mathbb C} K_{\mathbb C}$-decomposition 
for elements in $\Xi$. 

\par To start with let us introduce a {\it norm} on $G=SL(2,{\mathbb R})$ by 
setting 
$$\|g\|=\mbox{tr}(gg^t)\qquad (g\in G)\ . $$

\par Define a compact (here actually finite)
group 
$$L=N_{K_{\mathbb C}}(A_{\mathbb C}) \times Z= N_K(A)\times Z\ .$$
We let $L$ act on $A_{\mathbb C}$ via
$$(k,z)a=kak^{-1}z\qquad (k,z)\in L,\  a\in A_{\mathbb C}\ .$$
Then for every $x\in K_{\mathbb C} A_{\mathbb C} K_{\mathbb C}$ we 
have $x\in K_{\mathbb C} a_{K_{\mathbb C}}(x) K_{\mathbb C}$ with a unique continuous function 

$$a_{K_{\mathbb C}}: K_{\mathbb C} A_{\mathbb C} K_{\mathbb C} \to A_{\mathbb C}/ L, \ \  x\mapsto a_{K_{\mathbb C}}(x)\ . $$
In view of Theorem \ref{sph} and its proof we also obtain 
a holomorphic function 

\[ a_{K_{\mathbb C}}:\Xi \to A_{\mathbb C}/ L, \ \  x\mapsto 
\left\{ \begin{array}{rl}
a_{K_{\mathbb C}}(x) & \mbox{if}\  x\in K_{\mathbb C} A_{\mathbb C} K_{\mathbb C} / K_{\mathbb C},  \\
{\bf 1} & \mbox{if}\  x\not\in  K_{\mathbb C} A_{\mathbb C} K_{\mathbb C} / K_{\mathbb C}\ .\\
\end{array}\right. \]

Finally we define a {\it norm} on $A_{\mathbb C}$ by setting 
$$|\cdot|: A_{\mathbb C} \mapsto {\mathbb R}^+, \ \ a=\left(\begin{array}{cc} z & 0 \\ 0 & z^{-1} \end{array} \right)
\mapsto |z|^2 +|z|^{-2}\ . $$
Notice that $|\cdot |$ is $L$-invariant hence factors to 
a continuous positive function on $A_{\mathbb C}/ L$ also denoted by $|\cdot|$.

\begin{prop}\label{norm}
For all $g=h\exp(iX)K_{\mathbb C} \in \Xi$, ($h\in G, \ X\in \Omega$),  we have 
$$|a_{K_{\mathbb C}}(g)|\leq \|h\|\ .$$
\end{prop}

\begin{proof}
Using $G=KAK$ we can write 
\[\begin{array}{cccc}
g& = & \underbrace{\left( \begin{array}{cc} \cos r & \sin r \\ 
-\sin r & \cos r \end{array} \right)\left( \begin{array}{cc} 
s & 0 \\ 0 & s^{-1} \end{array} 
\right) \left( \begin{array}{cc} \cos t & \sin t \\ 
-\sin t & \cos t \end{array} \right) } & \underbrace{\left( 
\begin{array}{cc} 
e^{i\theta} & 0 \\ 0 & e^{-i\theta} \end{array} \right) } K_{\mathbb C} \\
& & =h \in KAK & \in \exp (i\Omega)
\end{array} \]
We may assume that $g\in K_{\mathbb C} A_{\mathbb C} K_{\mathbb C}/ K_{\mathbb C}$ and write 
$a_{K_{\mathbb C}}(g)=\left(\begin{array}{cc} z & 0 \\ 0 & z^{-1} \end{array} \right)$. 
Then 
\begin{eqnarray}
\Phi (g) & = & (s^{2}\cos^{2}t+s^{-2}\sin^{2}t)e^{2i\theta} + 
(s^{2}\sin^{2}t+s^{-2}\cos^{2}t)e^{-2i\theta} \nonumber \\
& = & \cos^{2}t(s^{2}e^{2i\theta}+s^{-2}e^{-2i\theta}) +
\sin^{2}t(s^{-2}e^{2i\theta}+s^{2}e^{-2i\theta}) \nonumber \\
& = & z^{2}+z^{-2} . 
\end{eqnarray}
Hence $\ z^{2}+z^{-2} \ $ is a convex combination of 
$\ s^{2}e^{2i\theta}+s^{-2}e^{-2i\theta} \ $ and
$\ s^{-2}e^{2i\theta}+s^{2}e^{-2i\theta} \ $  
and therefore lies on the vertical segment shown in Figure 2.
\\
\\

\epsfbox{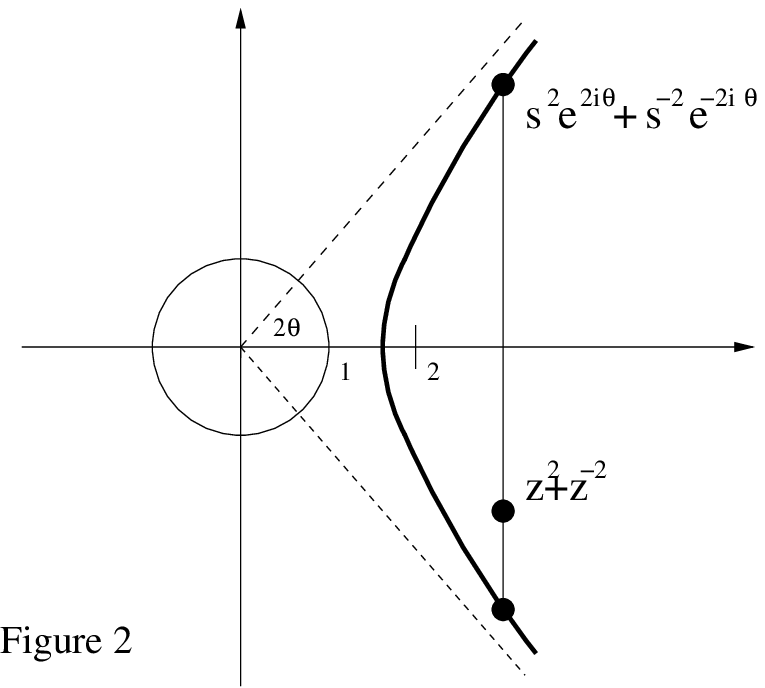}

The following Lemma \ref{normlem} applied to $u=s^2e^{2i\theta}$ and $v=z^2$ completes the proof
of the theorem. 
\end{proof}

\begin{lemma}\label{normlem} Let $r_1, r_2>0$ and $-{\pi\over 2}<\varphi_1, \varphi_2<{\pi\over 2}$. 
Define complex numbers by $u=r_1 e^{i\varphi_{1}}$ and $v=r_2 e^{i\varphi_{2}}$. 
Assume that 
$\ \Re(v+v^{-1}) = \Re(u+u^{-1}) \ $ and $\ \Im(v+v^{-1}) = \mu \Im(u+u^{-1})
 \ $ for some $\mu \in \mathbb{R}, \abs{\mu} \leq 1 $. 
Then $\quad \abs{\varphi_{2}} \leq \abs{\varphi_{1}} \ $ and $r_2+r_2^{-1}\leq r_1+r_1^{-1}$. 
\end{lemma}
\begin{proof}
The assumption written out in polar coordinates followed by squaring yields 
\begin{equation} \label{hirsch}
(r_{1}+r_{1}^{-1})^{2}\cos^{2}{\varphi_{1}}  =  (r_{2}+r_{2}^{-1})^{2}\cos^{2}{\varphi_{2}} \ ,
\end{equation}
\begin{equation*}
\mu^{2} (r_{1}-r_{1}^{-1})^{2}\sin^{2}{\varphi_{1}}  =  (r_{2}-r_{2}^{-1})^{2}\sin^{2}{\varphi_{2}} \ .
\end{equation*}

If $\abs{\varphi_{2}} > \abs{\varphi_{1}}$, then (\ref{hirsch}) implies
\[ (r_{1}+r_{1}^{-1})^{2} < (r_{2}+r_{2}^{-1})^{2} \ . \]
Since $\ (r_{1}+r_{1}^{-1})^{2} = (r_{1}-r_{1}^{-1})^{2}+4 \ $ this means
\[ (r_{1}-r_{1}^{-1})^{2} < (r_{2}-r_{2}^{-1})^{2} \ . \]
But now we end up with
\[ \mu^{2}(r_{1}-r_{1}^{-1})^{2}\sin^{2}{\varphi_{1}} < (r_{2}-r_{2}^{-1})^{2}\sin^{2}{\varphi_{2}} \ , \]
which contradicts the assumption. Thus $ \abs{\varphi_{2}} \leq \abs{\varphi_{1}}$ proving 
our first assertion. The second assertion is now immediate from (\ref{hirsch}).

\end{proof}

\bibliographystyle{amsplain}

\end{document}